\documentclass[11pt]{amsart} 
\usepackage[curve,v2]{xy}

\usepackage{amssymb,latexsym,amsfonts}

\newcommand{\tensor}{\otimes}

\renewcommand{\H}{\operatorname{Hilb}^2\hspace{-.05in}X}

\newcommand{\Tor}{\operatorname{Tor}}

\newcommand{\wtp}{\widetilde{\P}^n}

\newcommand{\wts}{\widetilde{\Sigma}}

\newcommand{\ses}[3]{0\rightarrow#1\rightarrow#2
   \rightarrow#3\rightarrow0}

\newcommand{\E}{{\mathcal E}}

\newcommand{\I}{{\mathcal I}}

\renewcommand{\O}{{\mathcal O}}
\renewcommand{\P}{{\mathbb{P}}}
\newcommand{\PP}{{\mathbb{P}}}

\newcommand{\N}{{\mathbb{N}}}
\newcommand{\NS}{{N_p^{\Sigma}}}

\renewenvironment{proof}{\par \medskip \noindent
{\sc Proof:}}{}

\begingroup
\newtheorem{thm}{Theorem}[section]   
\newtheorem{cor}[thm]{Corollary}     
\newtheorem{lemma}[thm]{Lemma}         
\newtheorem{prop}[thm]{Proposition}  

\newtheorem{conj}[thm]{Conjecture}        

\theoremstyle{definition}
\newtheorem{remark}[thm]{Remark}   
\newtheorem{notation}[thm]{Notation}
\newtheorem{notterm}[thm]{Notation and Terminology}

\endgroup

\newenvironment{rem}[2]{\refstepcounter{thm} \label{#2} 
\par \medskip \noindent {\bf #1 \thethm .}}{\par \medskip}

\begin{document}

\pagenumbering{arabic}

\title{Generation and syzygies of the first secant variety}




\author[Peter Vermeire]{Peter Vermeire}

\address{Department of Mathematics, 214 Pearce, Central Michigan
University, Mount Pleasant MI 48859}

\email{verme1pj@cmich.edu}



\date{\today}

\begin{abstract} 
Under certain effective positivity conditions, we show that the secant variety to a smooth variety 
satisfies $N_{3,p}$.  For smooth curves, we provide the best possible effective bound on the degree $d$ of the embedding, $d\geq 2g+3+p$.  
\end{abstract}

\maketitle

\section{Introduction}

\textbf{\textit{This should not be considered a final version.  Instead, I wanted to correct an error in version 2 of the posted preprint.  I was able to re-establish the degree bounds found there for curves, but I have not yet completed work on the higher-dimensional case.}}

We work throughout over an algebraically closed field of characteristic zero.
Let $X\subset\P^n$ be a smooth variety embedded by a line bundle $L$ and let $\Sigma_i$ denote the (complete) variety of $(i+1)$-secant $i$-planes.  Though secant varieties are a very classical subject, the majority of the work done involves determining the dimensions of secant varieties to well-known varieties.  Perhaps the two most well-known results in this direction are the solution by Alexander and Hirschowitz (completed in \cite{ah}) of the Waring problem for homogeneous polynomials and the classification of the Severi varieties by Zak \cite{zak}.  

More recently there has been great interest, e.g. related to algebraic statistics and algebraic complexity, in determining the equations defining secant varieties (e.g. \cite{ar}, \cite{bcg}, \cite{CGG05a}, \cite{CGG05b}, \cite{CGG07}, \cite{CGG}, \cite{unex}, \cite{cs}, \cite{gss}, \cite{kanev}, \cite{Lan06}, \cite{Lan08}, \cite{LM08}, \cite{LW07}, \cite{LW08}, \cite{sidsul}, \cite{ss}).  In this work, we use the detailed geometric information concerning secant varieties developed by Bertram \cite{bertram}, Thaddeus \cite{thaddeus}, and the author \cite{vermeireflip1} to lay some fundamental groundwork for studying not just the equations defining secant varieties, but the syzygies among those equations as well.

It was conjectured in \cite{eks} and it was shown in \cite{ravi} that if $C$ is a smooth curve embedded by a line bundle of degree at least $4g+2k+3$ then $\Sigma_k$ is set theoretically
defined by the $(k + 2)\times (k + 2)$ minors of a matrix of linear forms.  It was further shown in \cite{vermeireflip2} that if $X\subset\P^n$ satisfies condition $N_2$ then $\Sigma_1(v_d(X))$ is
set theoretically defined by cubics for $d\geq2$.

In \cite{vermeiresecantreg} it was shown that if $C$ is a smooth curve embedded by a line bundle of degree at least $2g+3$ then $\I_{\Sigma_1}$ is $5$-regular, and under the same hypothesis it was shown in \cite{sidver} that $\Sigma_1$ is arithmetically Cohen-Macaulay.  Together with the analogous well-known facts for the curve $C$ itself \cite{fujita}, \cite{mgreen}, \cite{mumford}, this led to the following conjecture, extending that found in \cite{vermeireflip2}:

\begin{conj}\cite{sidver}
Suppose that $C \subset \P^n$ is a smooth linearly normal curve of degree $d \geq 2g+2k+1+p$, where $p,k\geq0$.  Then 
\begin{enumerate}
\item $\Sigma_k$ is ACM and $\I_{\Sigma_k}$ has regularity $2k+3$ unless $g=0,$ in which case the regularity is $k+2$.  
\item $\beta_{n-2k-1, n+1} = \binom{g+k}{k+1}.$ 
\item $\Sigma_k$ satisfies $N_{k+2,p}$.\qed
\end{enumerate}
\end{conj}

\begin{remark}\label{whatsthatmean}
Recall \cite{EGHP} that a variety $Z\subset\P^n$ satisfies $\operatorname{N}_{r,p}$ if the ideal of $Z$ is generated in degree $r$ and the syzygies among the equations are linear for $p-1$ steps.  Note that condition $N_p$ \cite{mgreen} implies $N_{2,p}$.
\end{remark}

By the work of Green and Lazarsfeld \cite{mgreen},\cite{laz}, the conjecture holds for $k=0$.  Further, by \cite{fisher} and by \cite{vBH} it holds for $g\leq1$, and by \cite{sidver} parts (1) and (2) hold for $k=1$.  In this work, we show that part (3) holds for $k=1$ (Theorem~\ref{forcurves}).  More generally, we show that for an arbitrary smooth variety, $\Sigma_1$ satisfies $N_{3,p}$ for all sufficiently positive embeddings (Corollary~\ref{suffpos}), and we give effective results on 
arbitrary smooth varieties embedded by adjoint linear systems (Corollary~\ref{adjointsec}).  

Our approach combines the geometric knowledge of secant varieties mentioned above with the well-known Koszul approach of Green and Lazarsfeld.
To fix notation, if $L$ is a vector bundle on a projective variety $X$, then we let $\E_L=d_*(L\boxtimes\O)$, where $d:\operatorname{Bl}_{\Delta}(X\times X)\rightarrow \H$ is the natural double cover, and if $L$ is a globally generated line bundle on a projective variety $X$ inducing a morphism $f:X\rightarrow\P^n$, then we have the vector bundle $M_L=f^*\Omega_{\P^n}(1)$ on $X$.

\section{Preliminaries}

Our starting point is the familiar:

\begin{prop}\label{basic}
Let $X\subset\P^n$ be a smooth variety embedded by a line bundle $L$.  Then $\Sigma_1$ satisfies $\operatorname{N}_{3,p}$ if $H^1(\Sigma_1,\wedge^aM_L(b))=0$, $2\leq a\leq p+1$, $b\geq 2$.  

\end{prop}

\begin{proof}
Because $L$ also induces an embedding $\Sigma_1\subset\P^n$, we abuse notation and denote the associated vector bundle on $\Sigma_1$ by $M_L$.
Letting $F=\oplus\Gamma(\Sigma_1,\O_{\Sigma_1}(n))$ and
applying \cite[5.8]{eisenbud} to $\O_{\Sigma_1}$ gives the exact sequence:
\[
0 \to \Tor_{a-1}(F, k)_{a+b} \to H^1(\Sigma_1,\wedge^{a} M_L(b)) \to H^1(\Sigma_1,\wedge^{a}\O_{\P}^{r+1} \otimes \O_{\Sigma_1}(b))
\]
The vanishing in the hypothesis implies that $\Tor_1(F,k)_d=0$ for $d \geq k+1$, and hence that the first syzygies of $\O_{\Sigma_1}$, which are the generators of the ideal of $\Sigma_1$, are in degree $\leq k$.  The rest of the vanishings yield the analogous statements for higher syzygies.

\nopagebreak \hfill $\Box$ \par \medskip
\end{proof}

The technical portion of the paper is devoted to reinterpreting the vanishings in Proposition~\ref{basic} in terms of vanishings on the Hilbert scheme $\H$, and then finally on $X$ itself.

\begin{notterm}\label{maps}
Recall that an embedding $X\subset\P^n$ \textbf{\boldmath separates $k$ points} if every subscheme of $X$ of length $k$ spans a $\P^{k-1}\subset\P^n$.   A very ample line bundle $L$ is \textbf{\boldmath $k$-very ample} if the induced embedding separates $k+1$ points.  It is immediate that $k$-very ampleness implies $(k-1)$-very ampleness.

We will assume throughout that $X\subset \P^n$ is a $3$-very ample embedding of a smooth variety by $L=\O_X(1)$ that satisfies $N_{2,2}$.  For curves, an embedding given by a line bundle of degree at least $2g+3$ suffices \cite{mgreen}.  As we will be interested only in the first secant variety for the remainder of the paper, we write $\Sigma$ for $\Sigma_1$.

Under these hypotheses, the reader should keep in mind throughout the following morphisms \cite{vermeireflip1}
\begin{center}
{\begin{minipage}{1.5in}
\diagram
& & \H \\
& Z \cong \operatorname{Bl}_{\Delta}(X \times X) \urto^{d=\varphi|_Z} \dlto^{\pi_2} \dto^{\pi_1 = \pi|_Z} \ar@{^{(}->}[r]^{\hspace{.5in}i} & \wts \uto_{\varphi} \dto^{\pi}\\
X & X  \ar@{^{(}->}[r]  &\Sigma   
\enddiagram
\end{minipage}}
\end{center}
where 
\begin{itemize}
\item $\pi$ is the blow up of $\Sigma$ along $X$
\item  $i$ is the inclusion of the exceptional divisor of the blow-up
\item $d$ is the double cover, $\pi_i$ are the projections
\item $\varphi$ is the morphism induced by the linear system $|2H-E|$ which gives $\wts$ the structure of a $\P^1$-bundle over $\H$; note in particular that $\wts$ is smooth.  
\end{itemize}
We make frequent use of the rank $2$ vector bundle $\E_L=\varphi_*\O(H)=d_*\left(L\boxtimes\O\right)$, and note that $R^i\pi_*\O_{\wts}=H^i(X,\O_X)\tensor\O_X$ for $i\geq1$ (this is shown in \cite[Proposition 9]{vermeiresecantreg} for curves, but the same proof works in the general case).
\end{notterm}

\begin{prop}\label{prop: wts trans}
If $X$ is a smooth variety embedded by a $3$-very ample line bundle $L$ satisfying $N_{2,2}$, then $\Sigma$ satisfies $\operatorname{N}_{3,p}$ if $$H^1(\wts,\pi^*\wedge^aM_L(b))\rightarrow H^0(\Sigma,\wedge^aM_L(b)\tensor R^1\pi_*\O_{\wts})$$ is injective for $2\leq a\leq p+1$, $b\geq 2$.

\end{prop}

\begin{proof}
This follows immediately from the start of the 5-term sequence associated to the Leray-Serre spectral sequence:
$$0\rightarrow H^1(\Sigma,\wedge^aM_L(b))\rightarrow H^1(\wts,\pi^*\wedge^aM_L(b))\rightarrow H^0(\Sigma,\wedge^aM_L(b)\tensor R^1\pi_*\O_{\wts})$$
and Proposition~\ref{basic}.
\qed
\end{proof}


To simplify notation, we introduce a technical condition:

\begin{notation}\label{starstar}
For $p\geq 1$, we say a line bundle $L=\O_X(1)$ on $X\subset\P^n$ satisfies $\NS$ if
\begin{enumerate}
\item $L$ is $3$-very ample and satisfies $N_{2,p}$.
\item $H^i(\wts,\O_{\wts}(bH-E))=0$ for $i,b\geq1$.
\end{enumerate}
\end{notation}

As the vanishing condition in the definition of $\NS$ is non-trivial to understand, we explore several cases where it is satisfied in the next section.

\section{Condition $\NS$}

For curves, verification of $\NS$ is straightforward.
\begin{prop}\label{starstarforcurves}
Let $X\subset\P^n$ be a smooth curve satisfying $N_{p}$, $p\geq2$, with $L=\O_X(1)$ non-special.  Then $L$ satisfies $\NS$.
\end{prop}

\begin{proof}
We need to show $H^i(\wts,\O_{\wts}(bH-E))=0$ for $i,b\geq1$.

Because $X$ is projectively normal we have $H^i(\wtp,\O_{\wtp}(bH-E))=0$ for $i>0$, $b\geq1$.  Thus $H^i(\wts,\O_{\wts}(bH-E))=H^{i+1}(\wtp,\O_{\wtp}(bH-E)\tensor\I_{\wts})$.  By \cite[2.4(6)]{sidver}, we know that $H^{i+1}(\wtp,\O_{\wtp}(bH-E)\tensor\I_{\wts})=H^{i+1}(\P^n,\I_{\Sigma}(b))$.  

Now, for $i\geq1$, the arguments in \cite{vermeiresecantreg} and in \cite{sidver} go through under the stated hypotheses to give $H^{i+1}(\P^n,\I_{\Sigma}(b))=0$ for $b\geq1$.  The extra hypothesis used in those papers (namely, that $\deg(L)\geq 2g+3$) is needed only to show $H^{1}(\P^n,\I_{\Sigma}(b))=0$ for $b\geq1$.
\qed
\end{proof}

Verifying condition $\NS$ in the general case takes somewhat more work, but the end results are reasonable.
We first need a computation which will be used in both Proposition~\ref{vanishonsec} and in Theorem~\ref{forgeneral}.
\begin{lemma}\label{det}
Let $X$ be a smooth variety embedded by a $3$-very ample line bundle $L$ satisfying $N_{2,2}$.  Then $d^*\wedge^2\E_L=L\boxtimes L(-E_{\Delta})$.
\end{lemma}

\begin{proof}
Consider the sequence on $\wts$:
$$\ses{\O_{\wts}(-E)}{\O_{\wts}}{\O_Z}$$
As $R^0\varphi_*\O_{\wts}(-E)=0$, pushing down to $\H$ we have (\cite[3.10]{sidver})
$$\ses{\O_{\H}}{\O_{\H}\oplus M}{R^1\varphi_*\O_{\wts}(-E)}$$
where $d^*M=\O_{Z}(-E_{\Delta})$.

Thus $R^1\varphi_*\O_{\wts}(-E)=M$.  However, we know by \cite[5.1.2]{conrad} that $$\left(R^1\varphi_*\O_{\wts}(-E)\right)^*=R^0\varphi_*\left(\omega_{\wts/\H}\tensor\O_{\wts}(E)\right)$$ where $\omega_{\wts/\H}=\varphi^*\wedge^2\E_L(-2H)$ \cite[Ex.III.8.4b]{hart}.  Thus we have
$$M^*=\wedge^2\E_L\tensor\O_{\H}(-1)$$
and so $\varphi^*\wedge^2\E_L=\O_{\wts}(2H-E)\tensor\varphi^*M^*$.  Restricting (pulling back) this equality to $Z$ and noting (\cite[3.6]{vermeireidealreg}) that $\O_{Z}(2H-E)=L\boxtimes L(-2E_{\Delta})$, we have $d^*\wedge^2\E_L=L\boxtimes L(-E_{\Delta})$.
\qed
\end{proof}

We now interpret the vanishing condition in the definition of $\NS$ in terms of $X$.

\begin{prop}\label{vanishonsec}
Let $X\subset\P^n$ be a smooth variety embedded by a $3$-very ample line bundle $L$ satisfying $N_{2,2}$ such that $H^i(X\times X,L^{r+s}\boxtimes L^r\tensor \I_{\Delta}^{q})=0$ for $i,r\geq1$, $s\geq0$, $0\leq q\leq2r$.  Then $H^i(\wts,\O_{\wts}(bH-E))=0$ for $i,b\geq1$.
\end{prop}

\begin{proof}
Suppose $b=2r$ is even.  We know by the proof of \cite[3.6]{vermeireidealreg} that $\O_Z(bH-E)=L^{b-1}\boxtimes L\tensor\O(-2\Delta)$; thus
$$H^i(Z,\O_{Z}(bH-rE))=H^i(X\times X,L^r\boxtimes L^r\tensor \I_{\Delta}^{2r})=0$$  Because $\O_{\wts}(bH-rE)=\varphi^*\O_{\H}(r)$, we know $d_*\O_Z(bH-rE)=\O_{\H}(r)\tensor(\O\oplus M)$ for some line bundle $M$, and hence we know that $H^i(\H,\O_{\H}(r))=0$, but this says that 
$H^i(\wts,\O_{\wts}(bH-rE))=0$. From the sequences $$\ses{\O_{\wts}(bH-(k+1)E)}{\O_{\wts}(bH-kE)}{\O_{Z}(bH-kE)}$$
for $k+1\leq r$ we see that $H^i(\wts,\O_{\wts}(bH-E))=0$, as the cohomology of the rightmost terms vanishes by hypothesis since $H^i(Z,\O_{Z}(bH-kE))=H^i(X\times X,L^{b-k}\boxtimes L^k\tensor \I_{\Delta}^{2k})=0$.


Now, suppose that $b=2r+1$ is odd.  As in the previous paragraph, we have $\O_{\wts}((b-1)H-rE)=\varphi^*\O_{\H}(r)$, thus we see that $\varphi_*\O_{\wts}(bH-rE)=\O_{\H}(r)\tensor\varphi_*\O_{\wts}(H)=\O_{\H}(r)\tensor \E$.  It is therefore enough to show that $H^i(\H,\O_{\H}(r)\tensor \E)=0$, and then repeating the same argument as above gives $H^i(\wts,\O_{\wts}(bH-E))=0$.

We have the sequence on $Z$
$$\ses{K}{d^*\E_L}{L\boxtimes\O}$$
where $K=d^*\wedge^2\E_L\tensor \left(L^*\boxtimes\O\right)=\O\boxtimes L(-E_{\Delta})$ by Lemma~\ref{det}.  As in the proof of Lemma~\ref{det}, we have $d_*d^*\E_L=\E_L\oplus(\E_L\tensor M)$, thus $$d_*\left(\O_Z(2rH-rE)\tensor d^*\left(\E_L\tensor M^*\right)\right)=\E_L\tensor M^*(r)\oplus \E_L(r)$$
Thus it suffices to show $H^i(Z,\O_Z(2rH-rE)\tensor d^*\left(\E_L\tensor M^*\right))=0$.  However, we have 
$$K\tensor \O_Z(2rH-rE)\tensor d^*M^*=L^r\boxtimes L^{r+1}(-2rE_{\Delta})$$
and
$$L\boxtimes\O \tensor \O_Z(2rH-rE)\tensor d^*M^*=L^{r+1}\boxtimes L^r((-2r+1)E_{\Delta})$$
and so the cohomology of each vanishes by hypothesis.

\qed
\end{proof}

Fortunately, the vanishing in Proposition~\ref{vanishonsec} is not too difficult to understand.


\begin{prop}\label{uniform}
Let $X$ be a smooth variety of dimension $d$, $M$ a very ample line bundle.  Choose $k$ so that $k\geq d+3$ and so that $M^{k-d-1}\tensor\omega_X^*$ is big and nef.  Letting $L=M^k$, we have $$H^i(X\times X,L^{r+s}\boxtimes L^r\tensor \I_{\Delta}^{q})=0$$
for $i,r\geq1$, $s\geq0$, $0\leq q\leq2r$.  
\end{prop}

\begin{proof}
Note as above that $H^i(X\times X,L^{r+s}\boxtimes L^r\tensor \I_{\Delta}^{2r})=H^i(Z,L^{r+s}\boxtimes L^r\tensor \O(-2rE_{\Delta}))$, where $E_{\Delta}\rightarrow\Delta$ is the exceptional divisor of the blow-up.  Note further that $K_Z=K_X\boxtimes K_X\tensor \O((\dim X-1)E_{\Delta})$.



Assume first that $r\geq2$.
Then $$L^{r+s}\boxtimes L^r\tensor \O(-2rE_{\Delta})=K_Z\tensor (L^{r+s}-K_X)\boxtimes(L^r-K_X)\tensor \O((-d+1-q)E_{\Delta})$$
but this is $K_Z+B$ where
$$B=\left[(L-K_X)\boxtimes(L-K_X)\right]\tensor \left[L^{r+s-1}\boxtimes L^{r-1}\tensor\O\left((-d+1-q)E_{\Delta}\right)\right]$$

Because $M^k-K_X$ is ample, $(L-K_X)\boxtimes(L-K_X)$ is ample.  We are thus left to show that 
$$M^{k(r+s-1)}\boxtimes M^{k(r-1)}\tensor\O\left((-d+1-q)E_{\Delta}\right)$$
is globally generated.  However, as $k\geq d+3$, we have $k(r-1)\geq d-1+2r$ and so $M^{k(r+s-1)}\boxtimes M^{k(r-1)}\tensor\O\left((-d+1-2r)E_{\Delta}\right)$ is globally generated by \cite[3.1]{bel}.  Thus $B$ is big and nef and so vanishing follows from Kawamata-Viehweg vanishing \cite{kaw},\cite{vie}.

Now let $r=1$.  Then 
$$L^{1+s}\boxtimes L\tensor \O(-2E_{\Delta})=K_Z\tensor (M^{k+ks}-K_X)\boxtimes(M^{k}-K_X)\tensor \O((-d+1-q)E_{\Delta})$$
but this is $K_Z+B$ where
$$B=\left[(M^{k-d-1}-K_X)\boxtimes (M^{k-d-1}-K_X)\right]\tensor \left[M^{ks+d+1}\boxtimes M^{d+1}\right]\tensor \O((-d+1-q)E_{\Delta})$$
As above, $B$ is big and nef.
\qed
\end{proof}

\begin{remark}\label{adjoint}
There are numerous ways to rearrange the terms in Proposition~\ref{uniform} to produce the desired vanishing.  

For example, a similar argument shows that if $M$ is very ample, $\omega_X\tensor M$ is big and nef, and $B$ is nef, then letting $L=\omega_X\tensor M^k\tensor B$ gives the vanishing for $k\geq d+2$ (Cf. \cite[Theorem 1]{el}).
If, further, $B$ is also big, then letting $L=\omega_X\tensor M^k\tensor B$ gives the vanishing for $k\geq d+1$.
\end{remark}

\begin{remark}\label{fano}
In Proposition~\ref{uniform}, if $\omega_X^*$ is big and nef (e.g. $X$ is Fano) then a slight revision of the argument shows it is enough to take $L=M^k$ for $k\geq d+1$.
\end{remark}









\begin{remark}\label{gaussianetal}
Note that the vanishing condition in Proposition~\ref{vanishonsec} is intimately related to the surjectivity of the higher-order Gauss-Wahl maps as defined in \cite{wahl}.  
\end{remark}

\section{Main Results}


\begin{prop}\label{prop: wts van}
Let $X\subset\P^n$ be a smooth variety embedded by a line bundle $L$ satisfying $N_{p}^{\Sigma}$ with $H^i(X,L^k)=0$ for $i,k\geq1$.  Then $\Sigma$ satisfies $\operatorname{N}_{3,p}$ if $H^i(\wts,\pi^*\wedge^{a-1+i}M_L\tensor\O(2H-E))=0$ for $2\leq a\leq p+1$, $i\geq1$.


\end{prop}

\begin{proof}
We use Proposition~\ref{prop: wts trans}.  From the sequence on $\wts$
$$\ses{\pi^*\wedge^aM_L(bH-E)}{\pi^*\wedge^aM_L(bH)}{\pi^*\wedge^aM_L(bH)\tensor\O_Z}$$
we know \begin{eqnarray*}
H^1(Z,\pi^*\wedge^aM_L(bH)\tensor\O_{Z})&=&H^1\left(Z,\left(\wedge^aM_L\tensor L^b\right)\boxtimes \O_X\right)\\
&=&H^1\left(X\times X,\left(\wedge^aM_L\tensor L^b\right)\boxtimes \O_X\right)\\
&=&H^1(X,\O_X)\tensor H^0(X,\wedge^aM_L\tensor L^b).
\end{eqnarray*}
The first equality follows as the restriction of $\pi^*\wedge^aM_L(bH)$ to $Z$ is $\wedge^aM_L(bH) \boxtimes \O_X$, the second is standard, and for the third we use the K\"unneth formula together with the fact that $h^1(X, \wedge^aM_L \otimes L^b)=0$ as $X$ satisfies $N_{2,p}$. 

Thus $$h^1(\Sigma,\wedge^aM_L(b))=\operatorname{Rank}\left(H^1(\wts,\pi^*\wedge^aM_L(bH-E))\rightarrow H^1(\wts,\pi^*\wedge^aM_L(bH))\right)$$
and so by Proposition~\ref{prop: wts trans} it is enough to show that $H^1(\wts,\pi^*\wedge^aM_L\tensor\O(bH-E))=0$ for $2\leq a\leq p+1$, $b\geq2$.

From the sequence 
$$\ses{\pi^*\wedge^{a+1}M_L\tensor\O(bH-E)}{\wedge^{a+1}\Gamma\tensor\O(bH-E)}{\pi^*\wedge^aM_L\tensor\O((b+1)H-E)}$$
and the fact that $H^i(\wts,\O(bH-E))=0$, we see that $H^1(\wts,\pi^*\wedge^{a}M_L\tensor\O(bH-E))=H^{b-2}(\wts,\pi^*\wedge^{a+b-2}M_L\tensor\O(2H-E))$ for $b\geq2$.




\qed
\end{proof}

\begin{lemma}\label{downtohilb}
Let $X$ be a smooth variety embedded by a $3$-very ample line bundle $L$ satisfying $N_{2,2}$ and consider the morphism $\varphi:\wts\rightarrow \H\subset\P^s$ induced by the linear system $|2H-E|$.  Then $\varphi_*\wedge^aM_L=\wedge^aM_{\mathcal{E_L}}$, and hence $H^i(\wts,\pi^*\wedge^aM_L\tensor\O(2H-E))=H^i(\H,\wedge^aM_{\E_L}\tensor \O_{\H}(1))$.
\end{lemma}

\begin{proof}
Consider the diagram on $\wts$:
\begin{center}
{\begin{minipage}{1.5in}
\diagram
 &  &  & 0\dto & \\
 & 0\dto & 0\dto & K\dto & \\
0\rto &  \varphi^*M_{\E_L}\dto\rto & \Gamma(\H,\E_L)\otimes \O_{\wts}\dto\rto & \varphi^*\E_L\dto\rto & 0 \\
0\rto &  \pi^*M_{L}\dto\rto & \Gamma(X,L)\otimes \O_{\wts} \dto\rto & \pi^*L\dto\rto & 0 \\
 &  K\dto & 0 & 0 &  \\
 & 0 &  &  &
\enddiagram
\end{minipage}}
\end{center}
The vertical map in the middle is surjective as we have $\Gamma(\H, \E_L) = \Gamma(\wts,\O(H)) = \Gamma(X\times X,L\boxtimes\O) = \Gamma(X, L)$.  Therefore, surjectivity of the lower right horizontal map and commutativity of the diagram show that the righthand vertical map is surjective.  

Note that $R^i\varphi_*\varphi^*\E_L = \E_L \otimes R^i \varphi_*\O_{\wts}$ by the projection formula and that the higher direct image sheaves $R^i \varphi_*\O_{\wts}$ vanish as $\wts$ is a $\PP^1$-bundle over $\H.$ For the higher direct images, we have $R^i \varphi_*\pi^*L=0$ as the restriction of $L$ to a fiber of $\varphi$ is $\O(1)$ and hence the cohomology along the fibers vanishes.  From the rightmost column, we see $R^i\varphi_*K=0$. From the leftmost column, we have the sequence
$$\ses{\varphi^*\wedge^aM_{\E_L}}{\pi^*\wedge^aM_{L}}{\varphi^*\wedge^{a-1}M_{\E_L}\tensor K}$$
but as $R^i\varphi_*\left(K\tensor\varphi^*\wedge^{a-1}M_{\E_L}\right)=R^i\varphi_*K\tensor\wedge^{a-1}M_{\E_L}=0$, we have $\varphi_*\wedge^aM_L=\wedge^aM_{\mathcal{E}_L}$.
\qed
\end{proof}

Combining Proposition~\ref{prop: wts van} with Lemma~\ref{downtohilb} yields:

\begin{cor}\label{onhilb}
Let $X$ be a smooth variety embedded by a line bundle $L$ satisfying $N_{p}^{\Sigma}$ with $H^i(X,L^k)=0$ for $i,k\geq1$.
Then $\Sigma$ satisfies $\operatorname{N}_{3,p}$ if $$H^i(\H,\wedge^{a-1+i}M_{\E_L}\tensor\O(1))=0$$ for $2\leq a\leq p+1$, $i\geq1$.\qed
\end{cor}

\subsection{Curves}

We need a technical lemma, completely analogous to \cite[1.4.1]{laz}.

\begin{lemma}\label{belikerob}
Let $X\subset\P^n$ be a smooth curve embedded by a non-special line bundle $L$ satisfying $N_{2,2}$, let $x_1,\cdots,x_{n-2}$ be a general collection of distinct points, and let $D=x_1+\cdots+x_{n-2}$. Then there is an exact sequence of vector bundles on $X\times X$ $$\ses{L^{-1}(D)\boxtimes L^{-1}(D)(\Delta)}{d^*M_{\E_L}}{\displaystyle \bigoplus_i\left(\O(-x_i)\boxtimes\O(-x_i)\right)}$$
\end{lemma}

\begin{proof}
Choose a general point $x_1\in X$ and consider the following diagram on $X\times X$:
\begin{center}
{\begin{minipage}{1.5in}
\diagram
 & 0\dto & 0\dto & 0\dto & \\
0\rto & d^*M_{\E_{L(-x_1)}}\dto\rto & M_{L(-x_1)}\boxtimes \O\dto\rto & \left(\O\boxtimes L(-x_1)\right)(-\Delta)\dto\rto & 0\\
0\rto &  d^*M_{\E_{L}}\dto\rto & M_L\boxtimes\O \dto\rto & \O\boxtimes L(-\Delta)\dto\rto & 0 \\
0\rto &  \O(-x_1)\boxtimes\O(-x_1)\dto\rto & \O(-x_1)\boxtimes\O\dto\rto & \O\boxtimes\left(L\tensor\O_{x_1}\right) (-\Delta)\dto\rto & 0 \\
 &  0 & 0 & 0 & 
\enddiagram
\end{minipage}}
\end{center}
where the center column comes from \cite[1.4.1]{laz}.  Following just as in that proof, we obtain 
$$\ses{d^*M_{\E_{L(-D)}}}{d^*M_{\E_L}}{\displaystyle \bigoplus_i\left(\O(-x_i)\boxtimes\O(-x_i)\right)}$$
from the left column.  Note however that $d^*M_{\E_{L(-D)}}$ is a line bundle, and hence by Lemma~\ref{det} we see $d^*M_{\E_{L(-D)}}=\wedge^2\E^*_{L(-D)}=L^{-1}(D)\boxtimes L^{-1}(D)(\Delta)$.
\qed
\end{proof}

\begin{thm}\label{forcurves}
Let $X\subset\P^n$ be a smooth curve embedded by a line bundle $L$ with $\deg(L)\geq 2g+p+3$.  Then $\Sigma$ satisfies $\operatorname{N}_{3,p}$.
\end{thm}

\begin{proof}
We verify the condition in Corollary~\ref{onhilb}.
Pulling the sequence on $\H$
$$\ses{M_{\E_L}}{\Gamma(\H,\E_L)}{\E_L}$$
back to $Z=X\times X$ yields the diagram
\begin{center}
{\begin{minipage}{1.5in}
\diagram
 &  &  & 0\dto & \\
 & 0\dto & 0\dto & K\dto & \\
0\rto &  d^*M_{\E_L}\dto\rto & d^*\Gamma(\H, \E_L)\dto\rto & d^*\E_L\dto\rto & 0 \\
0\rto &  M_{L\boxtimes\O}\dto\rto & \Gamma(X,L)\dto\rto & L\boxtimes\O\dto\rto & 0 \\
 &  K\dto & 0 & 0 &  \\
 & 0 &  &  &
\enddiagram
\end{minipage}}
\end{center}
As in Lemma~\ref{det}, we have $d_*\O_{Z}=\O_{\H}\oplus M$ where $d^*M=\O(-\Delta)$, $d_*K=\E_L\tensor M$, and $K=d^*\wedge^2\E_L\tensor \left(L^*\boxtimes\O\right)=\O\boxtimes L(-\Delta)$.  From the left vertical sequence we have
$$\ses{\wedge^ad^*M_{\E_L}}{\wedge^aM_{L\boxtimes\O}}{\wedge^{a-1}d^*M_{\E_L}\tensor K}$$
and pushing down to $\H$ yields
$$\ses{\wedge^aM_{\E_L}\oplus \left(\wedge^aM_{\E_L}\tensor M\right)}{d_*\wedge^aM_{L\boxtimes\O}}{\wedge^{a-1}M_{\E_L}\tensor \E_L\tensor M}$$
Twisting this sequence by $\O_{\H}(1)\tensor M^*$ gives
$$\ses{\wedge^aM_{\E_L}(1)\tensor M^*\oplus \wedge^aM_{\E_L}(1)}{\O_{\H}(1)\tensor M^*\tensor d_*\wedge^aM_{L\boxtimes\O}}{\wedge^{a-1}M_{\E_L}(1)\tensor \E_L}$$
Since $d^*\O_{\H}(1)\tensor M^*=L\boxtimes L\tensor\O(-\Delta)$, it suffices (as in Proposition~\ref{vanishonsec}) to show that $$H^i(Z,\wedge^{a-1+i}d^*M_{\E_L}\tensor L\boxtimes L\tensor\O(-\Delta))=0$$ for $2\leq a\leq p+1$, $i=1,2$.

Now, by Lemma~\ref{belikerob} we have exact sequences
$$\ses{\wedge^{r-1}Q\tensor\O(D)\boxtimes\O(D)}{\wedge^rd^*M_{\E_L}\tensor L\boxtimes L\tensor\O(-\Delta)}{\wedge^rQ\tensor L\boxtimes L\tensor\O(-\Delta)}$$
where $Q=\displaystyle \bigoplus_i\left(\O(-x_i)\boxtimes\O(-x_i)\right)$.  

On the right, we have a direct sum of vector bundles of the form $F\boxtimes F(-\Delta)$ where $F$ is a line bundle of degree $\deg(L)-r$.  Thus $H^1$ and $H^2$ of the right side will vanish when $\deg(L)-r\geq 2g+1$.

On the left, we have a direct sum of vector bundles of the form $F\boxtimes F$ where $F$ is a line bundle of degree $n-2-(r-1)=\deg(L)-g-r-1$.  Because $x_1,\cdots,x_{n-2}$ are general, $H^1$ and $H^2$ of the left side will vanish when $\deg(L)-g-r-1\geq g$.  Combining these, we see that $H^i(Z,\wedge^{a-1+i}d^*M_{\E_L}\tensor L\boxtimes L\tensor\O(-E_{\Delta}))=0$ for $2\leq a\leq p+1$, $i=1,2$ as long as $\deg(L)\geq 2g+p+3$.

\qed
\end{proof}

\subsection{Higher Dimensions}

\begin{thm}\label{forgeneral}
Let $X\subset\P^n$ be a smooth variety embedded by a line bundle $L$ satisfying $\NS$, $p\geq1$, with $H^i(X,L^k)=0$ for $i,k\geq1$.
If $H^i(X,N^*_{X/\P^n}\tensor\wedge^{a-1+i}M_L\tensor L^2)=0$ for $2\leq a\leq p+1$ and for $i\geq1$, then $\Sigma$ satisfies $\operatorname{N}_{3,p}$. 
\end{thm}

\begin{proof}

As in Theorem~\ref{forcurves}, we have
$$\ses{\wedge^ad^*M_{\E_L}}{\wedge^aM_{L\boxtimes\O}}{\wedge^{a-1}d^*M_{\E_L}\tensor K}$$
and pushing down to $\H$ yields
$$\ses{\wedge^aM_{\E_L}\oplus \left(\wedge^aM_{\E_L}\tensor M\right)}{\wedge^aM_{\E_L}\oplus \left(F_a^2\tensor M\right)}{\wedge^{a-1}M_{\E_L}\tensor \E_L\tensor M}$$
where $F_a^2$ comes from the standard filtration
$$0\subset \wedge^aM_{\E_L}\subset F_a^2\subset\wedge^a\Gamma(\H, \E_L)$$
of $\wedge^a\Gamma$ associated to $\ses{M_{\E_L}}{\Gamma(\H, \E_L)}{\E_L}$ where
$$F_a^2/\wedge^aM_{\E_L}=\wedge^{a-1}M_{\E_L}\tensor\E_L; \wedge^a\Gamma/F_a^2=\wedge^{a-2}M_{\E_L}\tensor\wedge^2\E_L.$$
Twisting by $\O_{\H}(1)$, we see it is enough to show that $$H^i(Z,\wedge^{a-1+i}d^*M_{L}\tensor L\boxtimes L\tensor\O(-2E_{\Delta}))=0$$ for $2\leq a\leq p+1$, $i\geq 1$.  However, it is well-known that for $L$ very ample we have $H^i(Z,\wedge^{a-1+i}d^*M_{L}\tensor L\boxtimes L\tensor\O(-2E_{\Delta}))=H^i(X,N^*_{X/\P^n}\tensor\wedge^{a-1+i}M_L\tensor L^2)$.

\qed
\end{proof}




To verify the new vanishing condition we have:
\begin{prop}\label{verifynew}
Let $X$ be a smooth variety of dimension $d$, $M$ a very ample line bundle.  Choose $k\in\N$ so that $M^{k-ad-a-1}\tensor\omega_X^*$ is big and nef.  Letting $L=M^k$, we have $$H^i(X,N^*_{X/\P^n}\tensor\wedge^aM_L \tensor L^2)=0$$ for $i\geq1$.  
\end{prop}

\begin{proof}
Consider the product of $a+2$ factors $X\times X\times\cdots\times X$.  Then (cf. \cite{inamdar})
$$\left(\pi_1\right)_*L\boxtimes L\boxtimes\cdots\boxtimes L\tensor\I^2_{\Delta_{1,2}}\tensor\I_{\Delta_{1,3}}\tensor\cdots\tensor\I_{\Delta_{1,a+2}}=M_L^{\tensor a}\tensor N^*(2)$$
Arguing as in Proposition~\ref{uniform}, we obtain $H^1(X,M_L^{\tensor a}\tensor N^*(2))=0$.  However, as we are working in characteristic $0$, $\wedge^aM_L \tensor N^*(2)$ is a summand of $M_L^{\tensor a}\tensor N^*(2)$.
\qed
\end{proof}




\begin{cor}\label{suffpos}
Let $X$ be a smooth variety, $M$ an ample line bundle, and embed $X$ by $M^k$.  Then for all $k>>0$, $\Sigma$ satisfies $N_{3,p}$.  
\end{cor}

\begin{proof}
Letting $L=M^k$ for $k>>0$, we know \cite{mgreen},\cite{inamdar} $L$ satisfies $N_{p}$. By Proposition~\ref{vanishonsec}, $L$ satisfies $\NS$.  Finally, by Proposition~\ref{verifynew}, we have $H^i(X,N^*_{X/\P^n}\tensor\wedge^{a-1+i}M_L\tensor L^2)=0$ for $2\leq a\leq p+1$ and for $i\geq1$.
\qed
\end{proof}




\begin{cor}\label{adjointsec}
Let $X\neq\P^d$ be a smooth projective variety of dimension $d$, $M$ a very ample line bundle such that $K_X\tensor M$ is ample.  Embedding $X$ by $L=K_X\tensor M^{(p+2)d+1}$, $p\geq1$, we have $\Sigma$ satisfies $N_{3,p}$.
\end{cor}

\begin{proof}
In \cite[3.1]{el} it is shown that $L=K_X\tensor M^{d+p+2}$ satisfies $N_{p+2}$.  The result now follows as in Remark~\ref{adjoint}.
\qed
\end{proof}

We conjecture what we believe to be the best possible result in general:
\begin{conj}
Let $X^d\subset\P^n$ be a smooth projective variety with $H^i(X,\O(k))=0$ for $i,k\geq1$ satisfying $N^{\Sigma}_{p+d+1}$.  Then $\Sigma$ satisfies $N_{3,p}$.
\qed
\end{conj}

\begin{rem}{Remark}{partial}
We can show that under the hypotheses of the Conjecture that $\Sigma$ satisfies $N_{3+d,p}$.
\qed
\end{rem}

\section{Acknowledgments}
This project grew out of work done together with Jessica Sidman, and has benefited greatly from her insight and input, as well as from her comments regarding a preliminary draft of this work.  I would also like to thank Lisa DeMeyer, Hal Schenck, Greg Smith, and Jonathan Wahl for helpful discussion and comments.

\end{document}